\documentclass[final]{amsart}
\usepackage{amssymb}
\usepackage{hhline}
\usepackage{graphicx}
\usepackage{stmaryrd}

\newtheorem{theorem}{Theorem}[section]

\newtheorem{proposition}[theorem]{Proposition}
\newtheorem{corollary}[theorem]{Corollary}

\newtheorem{problem}[theorem]{Problem}
\theoremstyle{definition}

\newcommand{\Aut}{\mathrm{Aut\mkern 2mu}}

\newcommand{\M}{\mathrm{M\mkern 1mu}}

\title{Classifications of dimonoids with at most three elements}
\author{Volodymyr M. Gavrylkiv}
\address[V.~Gavrylkiv]{Vasyl Stefanyk Carpathian National University, Ivano-Frankivsk, Ukraine} \email{vgavrylkiv@gmail.com}
\begin{document}

\begin{abstract}
In this paper, we present complete classifications, up to isomorphism, of all two-element dimonoids, all commutative three-element dimonoids, and all abelian three-element dimonoids. We show that, up to isomorphism, there exist exactly 8  two-element dimonoids, of which 3 are commutative. Among these, 4 are abelian, and the remaining nonabelian dimonoids form 2 pairs of dual dimonoids. Furthermore, there are exactly 5 pairwise nonisomorphic trivial dimonoids of order 2. For dimonoids of order 3, we prove that there are precisely 14 pairwise nonisomorphic commutative dimonoids, including 12 trivial dimonoids and a single pair of  nonabelian nontrivial dual dimonoids. We also establish that, up to isomorphism, there are 17 abelian dimonoids of order 3, consisting of 12 trivial commutative dimonoids and 5 noncommutative nontrivial ones. In addition, we demonstrate the existence of at least 26 pairwise nonisomorphic nonabelian noncommutative dimonoids of order 3. Among them, there are exactly 6 pairs of trivial dual dimonoids and at least 7 pairs of nontrivial dual dimonoids.
\end{abstract}

\maketitle

\section*{Introduction}

The notions of a dialgebra and a dimonoid were introduced by J.-L. Loday~\cite{Lod}. A {\em dimonoid} is an algebraic structure $(D,\dashv,\vdash)$ consisting of a  set $D$ equipped with two associative binary operations $\dashv$ and $\vdash$ satisfying the following axioms:
\begin{align*}
(x \dashv y) \dashv z = x \dashv (y \vdash z), \hspace{10mm}(D_1) \\
(x \vdash y) \dashv z = x \vdash (y \dashv z), \hspace{10mm} (D_2)\\
(x \dashv y) \vdash z = x \vdash (y \vdash z). \hspace{10mm} (D_3)
\end{align*}

Each semigroup $(D,\dashv)$ can naturally be regarded as a dimonoid $(D,\dashv,\dashv)$, referred to as the {\em trivial} dimonoid and denoted simply by $D$. In this way, dimonoids generalize semigroups. Dimonoids also admit linear analogues known as dialgebras. A {\em dialgebra} is a vector space over a field equipped with two  bilinear associative binary operations satisfying the axioms of a dimonoid. Consequently, many results concerning dimonoids have direct applications in dialgebra theory~\cite{B,F, Lod, M, ZCY}. In recent years, dimonoids have become standard tools in the study of various structures, particularly in the theory of Leibniz algebras. Notably, T.~Pirashvili~\cite{P} introduced the concept of a duplex, a generalization of dimonoids, and constructed the free duplex. The properties of free dimonoids were employed in~\cite{Lod} to characterize free dialgebras and to study their cohomologies. In~\cite{Liu}, the notion of a dimonoid was used to define and investigate one-sided dirings. Furthermore, dimonoids are closely related to restrictive bisemigroups~\cite{Sh} and doppelsemigroups~\cite{GR1, GDS2, GUDS, GSDS, GLDS, ZhDADM2016}.

One of the earliest foundational results on dimonoids is due to Loday~\cite{Lod}, who provided a description of the absolutely free dimonoid generated by a given set.
A wide range of classes of dimonoids have been systematically investigated by Anatolii Zhuchok and Yurii Zhuchok. In~\cite{ZhAL2011}, the independence of the dimonoid axioms was established. Commutative, free commutative, and free abelian dimonoids were studied in~\cite{ZhADM2009}, \cite{ZhADM2010}, and~\cite{ZhADM2015}, respectively. The structure of dibands of subdimonoids and semilattice decompositions of dimonoids was explored in~\cite{ZhMS2011, ZhDM2011}. Free rectangular dimonoids, as well as free normal and free $(lr, rr)$-dibands, were constructed in~\cite{ZhADM2011a}, \cite{ZhADM2011}, and~\cite{ZhADM2013}, respectively. Free abelian dibands and some of their properties were studied in~\cite{ZhVLU2017, ZhADM2018}. The least semilattice congruence on free dimonoids was described in~\cite{ZhUMJ2011}. Free products of dimonoids and relatively free dimonoids were the focus of several works, including~\cite{ZhQRL2013, ZhMP2014, ZhCA2017a, ZhIJAC2021}. Moreover, the free left $n$-nilpotent and free left $n$-dinilpotent dimonoids were constructed in~\cite{ZhADM2013nil, ZhSF2016}. Representations of ordered dimonoids via binary relations were examined in~\cite{ZhAEJM2014}.
Significant contributions to the theory of endomorphisms and automorphisms in the context of dimonoids were made by Y. Zhuchok in~\cite{ZhBAS2014,ZhCA2017,ZhJA2024}.

In~\cite{Gdim1}, we studied the properties of dual dimonoids and, within the class of noncommutative dimonoids, constructed various examples of abelian, nonabelian, and rectangular dimonoids. The algebraic structure of these dimonoids was examined in detail, including computations of their automorphism groups. In the present work, we build upon these results, as well as those obtained in~\cite{GR1, GDS2}, to provide a complete classification, up to isomorphism, of all two-element dimonoids, all commutative three-element dimonoids, and all abelian three-element dimonoids.

\section{Preliminaries on semigroups}

An element $e$ of a semigroup $(S,*)$ is called  {\em a left identity} (resp.  {\em a right identity}) in $S$ if $e*a=a$ (resp. $a*e=a$) for any $a\in S$. An element $e$  is called  {\em an identity} if $e$ is a left identity and a right identity.

Let $(S,*)$ be a semigroup and $e\notin S$. The binary operation $*$ defined on  $S$  can be extended to $S\cup\{e\}$ putting $e*s=s*e=s$ for all $s\in S\cup \{e\}$.
The notation $(S,*)^{+1}$  denotes a monoid $(S\cup\{e\}, *)$ obtained from $(S,*)$ by adjoining the extra identity $e$ (regardless of whether $(S,*)$ is or is not a monoid).

Let $(M, *)$ be a monoid with  identity $e$ and $M^{\tilde{1}}=M\cup\{\tilde{1}\}$, where $\tilde{1}\notin M$.
The binary operation $*$ defined on  $M$  can be extended to $M^{\tilde{1}}$
putting $\tilde{1}* m=m* \tilde{1}=m$ for all $m\in M$ and ${\tilde{1}}*{\tilde{1}}=e$.
The notation $(M,*)^{\tilde{1}}$  denotes the semigroup  obtained from $(M,*)$ by adjoining an extra element $\tilde{1}$. Note that  $(M,*)^{\tilde{1}}$ is not a monoid.

An element $e$ of a semigroup $(S,*)$ is called an {\em idempotent} if $e*e=e$.  The semigroup is a {\em band}, if all its elements are idempotents. Commutative bands are called {\em semilattices}.
By $L_n$ we denote the {\em linear semilattice} $\{0, 1,\ldots, n-1\}$ of order $n$, endowed with the operation of minimum.

\smallskip

A semigroup $(S,*)$ is called {\em monogenic} if it is generated by some element $a\in S$ in the sense that $S=\{a^n\}_{n\in\mathbb N}$. If a monogenic semigroup is
infinite then it is isomorphic to the additive semigroup $\mathbb N$ of positive integer numbers. A
finite monogenic semigroup $S=\langle a\rangle$ also has simple
structure, see  \cite{Howie}. There are positive integer numbers $r$
and $m$ called the {\em index} and the {\em period} of $S$ such
that
\begin{itemize}
    \item $S=\{a,a^2,\dots,a^{r+m-1}\}$ and $r+m-1=|S|$;
    \item $a^{r+m}=a^{r}$;
    \item $C_m:=\{a^r,a^{r+1},\dots,a^{r+m-1}\}$ is
    a cyclic and maximal subgroup of $S$ with the
    identity $e=a^n\in C_m$ and generator $a^{n+1}$, where $n\in (m\cdot\mathbb N)\cap\{r,\dots,r+m-1\}$.
\end{itemize}

We denote by $\M_{r,m}$ a finite monogenic semigroup of index $r$ and period $m$.

\smallskip

An element $z$ of a semigroup $S$ is called  {\em a left zero} (resp.  {\em a right zero}) in $S$ if $z*a=z$ (resp. $a*z=z$) for any $a\in S$. An element $0$  is called  {\em a  zero} if $0$ is a left zero and a right zero.

Let $(S,*)$ be a semigroup and $0\notin S$. The binary operation $*$ defined on  $S$  can be extended to $S\cup\{0\}$ putting $0*s=s*0=0$ for all $s\in S\cup \{0\}$.
The notation $(S,*)^{+0}$  denotes a semigroup $(S\cup\{0\},*)$ obtained from $(S,*)$ by adjoining the extra zero $0$ (regardless of whether $(S,*)$ has or has not a zero). 

\smallskip

A semigroup $(S,*)$ is called a {\em null semigroup} if there exists an element $0\in S$ such that $x*y=0$ for all $x,y\in S$. In this case  $0$ is a zero of $S$. All null semigroups on the same set are isomorphic. By $O_S$ we denote a null semigroup on a set $S$. If $S$ is finite of cardinality $|S|=n$, then instead of $O_S$ we use $O_n$. 

\smallskip

Let $S$ be a nonempty set, $0\in S$ and $A\subset S\setminus\{0\}$. Define  the binary operation $*$ on  $S$ in the following way:

$$x* y=\begin{cases}
x, \text{ if }y=x\in A \\
0,\ \text{otherwise}.
\end{cases}$$

It is easy to check that a set $S$ endowed with the operation $*$ is a commutative semigroup with  zero $0$, and we denote this semigroup by $O^A_S$. If $A=S\setminus\{0\}$, then $O^A_S$ is a semilattice. In the case when $A$ is an emptyset,  $O^A_S$ coincides with a null  semigroup  with  zero $0$. The semigroups $O^A_S$ and $O^B_T$ are isomorphic if and only if $|S|=|T|$ and $|A|=|B|$. If $S$ is a finite set of cardinality $|S|=n$ and $|A|=m$, then we use $O_n^m$ instead of $O^A_S$.

\smallskip

If $(S,*)$ is a semigroup, then the semigroup $(S,{*}^d)$ with operation $x{*}^d y=y* x$ is called {\em dual} to $(S,*)$, denoted $(S,*)^d$. It follows that $(S,*)^d = (S,*)$ if and only if $(S,*)$ is a commutative semigroup.

\smallskip

A semigroup $(S,*)$ is said to be a {\em left} (resp. {\em right}) {\em zero semigroup} if $a*b=a$ (resp. $a*b=b$) for any $a,b\in S$. By  $LO_S$ and $RO_S$ we denote  a left zero semigroup and a right zero semigroup on a set $S$, respectively.  It is easy to see that the semigroups $LO_S$ and $RO_S$ are dual. If $S$ is finite of cardinality $|S|=n$, then instead of $LO_S$ and $RO_S$ we use $LO_n$ and $RO_n$, respectively.

\smallskip

Let $S$ be a nonempty set, $A\subset S$ and $0\notin S$. Define  the binary operation $*$ on  $S^0=S\cup\{0\}$ in the following way:

$$x* y=\begin{cases}
x, \text{ if }y \in A \\
0, \text{ if }y \in S^0\setminus A.
\end{cases}$$

It is easy to check that a set $S^0$ endowed with the operation $*$ is a semigroup with  zero $0$, and we denote this semigroup by $LO^{\sim 0}_{A\leftarrow S}$. If $A=S$, then $LO^{\sim 0}_{A\leftarrow S}$ coincides with $LO^{+0}_S$. In the case when $A$ is an emptyset,  $LO^{\sim 0}_{A\leftarrow S}$ coincides with a null semigroup $O_{S^0}$ with  zero $0$. The semigroups $LO^{\sim 0}_{A\leftarrow S}$ and $LO^{\sim 0}_{B\leftarrow T}$ are isomorphic if and only if $|S|=|T|$ and $|A|=|B|$. If $S$ is a finite set of cardinality $|S|=n$ and $|A|=m$, then we use $LO^{\sim 0}_{m\leftarrow n}$ instead of $LO^{\sim 0}_{A\leftarrow S}$.

By $RO^{\sim 0}_{A\leftarrow S}$ we denote a dual semigroup of $LO^{\sim 0}_{A\leftarrow S}$.
 
\smallskip 

Let  $a$ and $c$ be different elements of a set $S$. Define the associative binary operation $\dashv_c^a$ on  $S$ in the following way:

$$x\dashv_c^a y=\begin{cases}
a,\text{ if } x=y=a \\
c,\text{ if } x=a\text{ and } y\neq a\\
x,\text{ if } x\neq a.
\end{cases}$$

If $|S|\geq 3$, then $(S,\dashv_c^a)$ is a noncommutative band in which all elements $z\neq a$ are left zeros.

It is not difficult to check that for any different $b, d\in S$, the semigroups $(S,\dashv_c^a)$ and $(S,\dashv_d^b)$ are
isomorphic. We denote this semigroup by $LOB_S$. If $S$ is a finite set of cardinality $|S|=n$,  then we use $LOB_n$ instead of $LOB_S$. 

By $ROB_S$ we denote a dual semigroup of $LOB_S$.

\smallskip

Let $S$ be a nonempty set, $A$ be a nonempty subset of $S$, and $a\in A$. Define the associative binary operation $*$ on  $S$ in the following way:

$$x* y=\begin{cases}
x,\text{ if } x\in A \\
a,\text{ if } x \notin A.
\end{cases}$$

We denote the semigroup $(S,*)$  by $LO_{A\leftarrow S}$. It follows that  all elements $z\in A$ are left zeros of $LO_{A\leftarrow S}$. If $A=\{a\}$, then $LO_{A\leftarrow S}$ coincides with a null semigroup $O_S$ with zero $a$. If $A=S$, then $LO_{A\leftarrow S}$ coincides with a left zero semigroup $LO_S$. The semigroups $LO_{A\leftarrow S}$ and $LO_{B\leftarrow T}$ are isomorphic if and only if $|S|=|T|$ and $|A|=|B|$. If $S$ is a finite set of cardinality $|S|=n$ and $|A|=m$, then we use $LO_{m\leftarrow n}$ instead of $LO_{A\leftarrow S}$.

By $RO_{A\leftarrow S}$ we denote a dual semigroup of $LO_{A\leftarrow S}$.

\bigskip

Following the algebraic tradition, we take for a model of the class of cyclic groups of order $n$ the multiplicative group
$C_n=\{z\in\mathbb C:z^n=1\}$ of $n$-th roots of $1$. For a set $X$ by $S_X$ we denote the group of all bijections of $X$.

\section{Some definitions and basic properties of dimonoids}

In this section, we recall several useful results on dimonoids and establish auxiliary propositions that will be frequently used in the subsequent investigations.

An element $0\in D$  is called  a {\em   zero of a dimonoid} $(D,\dashv, \vdash)$~\cite{ZhADM2013nil} if $0$ is a zero of $(D,\dashv)$ and a  zero of $(D,\vdash)$. Let $(D,\dashv, \vdash)$ be a dimonoid and $0\notin D$. The binary operations defined on  $D$  can be extended to $D\cup\{0\}$ putting $0\dashv d=d\dashv 0=0=0\vdash d=d\vdash 0 $ for all $d\in D\cup \{0\}$. The notation $(D,\dashv, \vdash)^{+0}$  denotes a dimonoid $D\cup\{0\}$ obtained from $D$ by adjoining the extra zero $0$. 
%It follows  that $\Halo((D,\dashv, \vdash)^{+0}) = \Halo(D,\dashv, \vdash)$.

\smallskip

A dimonoid $(D,\dashv, \vdash)$ is called {\em abelian}~\cite{ZhADM2015} if $x \dashv y = y \vdash x$ for all $x,y\in D$.

\smallskip

Let $(D,\dashv, \vdash)$ be a dimonoid. Define new operations $\dashv^d$ and  $\vdash^d$ on $D$ by 
$$x \dashv^d y = y \vdash x\ \ \text{  and  }\ \  x \vdash^d y = y \dashv x.$$
It is immediate to check that $(D,\dashv^d, \vdash^d)$ is a new dimonoid, called the {\em  dual dimonoid of $(D,\dashv, \vdash)$}~\cite{Lod}, which we denote by $(D,\dashv, \vdash)^d$. It follows that the unary duality operation is involutive in the sense that $((D,\dashv, \vdash)^d)^d=(D,\dashv, \vdash)$. In fact, $(D,\dashv, \vdash)^d$ is a dimonoid if and only if $(D,\dashv, \vdash)$ is a dimonoid.  As usual, a dimonoid $(D,\dashv, \vdash)$ is said to be {\em self-dual} if $(D,\dashv, \vdash)^d=(D,\dashv, \vdash)$. As established in~\cite{Gdim1}, a dimonoid $(D,\dashv, \vdash)$ is abelian if and only if it is self-dual, which in turn holds if and only if the semigroups $(D, \vdash)$ and $(D,\vdash)$ are dual to each other. Consequently, nonabelian dimonoids are divided into the pairs of dual dimonoids.

\smallskip

A dimonoid $(D,\dashv,\vdash)$ is called {\em commutative}~\cite{ZhADM2009} if both semigroups $(D,\dashv)$ and $(D,\vdash)$ are commutative.

\smallskip

Since commutative semigroups $(D,\dashv)$ and $(D,\vdash)$ are dual if and only if their operations coincide, all commutative nontrivial dimonoids are nonabelian.
On the other hand, it is clear to see that all commutative trivial dimonoids are abelian and all noncommutative trivial dimonoids are nonabelian. A left zero and a right zero dimonoid $(D,\dashv,\vdash)$ with operations $x\dashv y = x$ and $x\vdash y = y$ \cite{Lod} is an example of a nontrivial abelian noncommutative dimonoid. In the section~\ref{3eldim} we give examples of commutative nonabelian dimonoids, see also~\cite{ZhADM2009}.

\bigskip

The axioms $(D_1)$ and $(D_3)$ of a dimonoid imply the following proposition.

\begin{proposition}\label{rl_id} Let $(D,\dashv, \vdash)$ be a dimonoid. If a semigroup $(D,\dashv)$ contains a left identity or a semigroup $(D,\vdash)$ contains a right identity, then the operations of a dimonoid $(D,\dashv, \vdash)$ coincide.
\end{proposition}

\begin{proposition}\label{lz} Let $(D,\dashv, \vdash)$ be a dimonoid. If $z\in D$ is a left zero of a semigroup $(D,\vdash)$, then $z$ is a left zero of a semigroup $(D,\dashv)$ as well.
\end{proposition}

\begin{proof} Taking into account that for any $a\in D$ the following equalities hold
$$z \dashv a = (z \vdash a) \dashv a = z \vdash (a \dashv a) = z,$$
we conclude that $z$ is a left zero of a semigroup $(D,\dashv)$.
\end{proof}

Dually, we prove the following proposition.

\begin{proposition}\label{rz} Let $(D,\dashv, \vdash)$ be a dimonoid. If $z\in D$ is a right zero of a semigroup $(D,\dashv)$, then $z$ is a right zero of a semigroup $(D,\vdash)$ as well.
\end{proposition}

\begin{proof} Since for any $a\in D$ the following equalities hold
$$a \vdash z = a \vdash (a \dashv z) = (a \vdash a) \dashv z = z,$$
we conclude that $z$ is a right zero of a semigroup $(D,\vdash)$.
\end{proof}

\begin{corollary}\label{zerodim} Let $(D,\dashv, \vdash)$ be a commutative dimonoid. An element $z\in D$ is a zero of a semigroup $(D,\dashv)$ if and only if $z$ is a zero of a semigroup $(D,\vdash)$.
\end{corollary}

Propositions~\ref{lz} and~\ref{rz} imply the following corollary.

\begin{corollary}\label{rightcoin} Let $(D,\dashv, \vdash)$ be a dimonoid. If $(D,\dashv)$ is a right zero semigroup or $(D,\vdash)$ is a left zero semigroup, then the operations of a dimonoid $(D,\dashv, \vdash)$ coincide.
\end{corollary}

A bijective map $\psi : D_1 \to D_2$ is called an {\em isomorphism } from a dimonoid $(D_1,\dashv_1, \vdash_1)$ to a dimonoid $(D_2,\dashv_2, \vdash_2)$ if $$\psi(a\dashv_1 b)=\psi(a)\dashv_2\psi(b)\ \ \text{  and  }\ \ \psi(a\vdash_1 b)=\psi(a)\vdash_2\psi(b)$$ for all $a,b\in D_1$.

If there exists an isomorphism from a dimonoid $(D_1,\dashv_1, \vdash_1)$ to a dimonoid $(D_2,\dashv_2, \vdash_2)$, then $(D_1, \dashv_1, \vdash_1)$ and $(D_2, \dashv_2, \vdash_2)$ are said to be {\em isomorphic}, denoted $(D_1,\dashv_1, \vdash_1)\cong (D_2,\dashv_2, \vdash_2)$. An isomorphism $\psi: D\to D$ is called an {\em   automorphism} of a dimonoid $(D,\dashv, \vdash)$. By $\Aut(D,\dashv, \vdash)$ we denote the automorphism group of a dimonoid $(D,\dashv, \vdash)$. It follows that $\Aut((D,\dashv, \vdash)^{+0}) = \Aut(D,\dashv, \vdash)$.

For a dimonoid $(D,\dashv,\vdash)$, if $\mathbb S$ and $\mathbb T$ denote the semigroups $(D,\dashv)$ and $(D,\vdash)$, respectively, then $\mathbb S \rbag \mathbb T$ stands for the dimonoid $(D,\dashv,\vdash)$.

%Let $(D_1,\dashv_1, \vdash_1)$ be such a dimonoid that for each dimonoid $(D_2,\dashv_2, \vdash_2)$ the isomorphisms 
%$(D_2,\dashv_2)\cong (D_1,\dashv_1)$ and $(D_2,\vdash_2)\cong (D_1,\vdash_1)$ imply $(D_2,\dashv_2, \vdash_2)\cong
%(D_1,\dashv_1, \vdash_1)$. If $\mathbb S$ and $\mathbb T$ are model semigroups of classes of semigroups isomorphic to
%$(D_1,\dashv_1)$ and $(D_1, \vdash_1)$, respectively, then by $\mathbb S\rbag\mathbb T$  we denote a model dimonoid of the class of dimonoids isomorphic to $(D_1,\dashv_1, \vdash_1)$.

\begin{proposition}\label{isoLO} Let $(D_1,\dashv_1, \vdash_1)$ and $(D_2,\dashv_2, \vdash_2)$ be  dimonoids such that  and $(D_1,\dashv_1)$ and $(D_2,\dashv_2)$ are left zero semigroups.  Dimonoids $(D_1,\dashv_1, \vdash_1)$ and $(D_2,\dashv_2, \vdash_2)$ are isomorphic if and only if semigroups $(D_1,\vdash_1)$ and $(D_2,\vdash_2)$ are isomorphic.
\end{proposition}

\begin{proof} It is immediate to observe that if dimonoids $(D_1,\dashv_1, \vdash_1)$ and $(D_2,\dashv_2, \vdash_2)$ are isomorphic, then semigroups $(D_1,\vdash_1)$ and $(D_2,\vdash_2)$ are also isomorphic.
Conversely, let $\psi: D_1 \to D_2$ be an isomorphism from a semigroup $(D_1,\vdash_1)$ to a semigroup $(D_2,\vdash_2)$. Then, necessarily, $|D_1| = |D_2|$. Taking into account that any bijective map is an isomorphism from a left zero semigroup $(D_1,\dashv_1)$ to a left zero semigroup $(D_2,\dashv_2)$, it follows that $\psi$ is also an isomorphism from a left zero semigroup $(D_1,\dashv_1)$ to a left zero semigroup $(D_2,\dashv_2)$. Therefore, $\psi$ is an isomorphism from a dimonoid $(D_1,\dashv_1, \vdash_1)$ to a dimonoid $(D_2,\dashv_2, \vdash_2)$.
\end{proof}

Dually, one can prove the following proposition.

\begin{proposition}\label{isoRO} Let $(D_1,\dashv_1, \vdash_1)$ and $(D_2,\dashv_2, \vdash_2)$ be  dimonoids suct that  and $(D_1,\vdash_1)$ and $(D_2,\vdash_2)$ are right zero semigroups.  Dimonoids $(D_1,\dashv_1, \vdash_1)$ and $(D_2,\dashv_2, \vdash_2)$ are isomorphic if and only if semigroups $(D_1,\dashv_1)$ and $(D_2,\dashv_2)$ are isomorphic.
\end{proposition}

\begin{proposition}\label{isoabelian} Let $(D_1,\dashv_1, \vdash_1)$ and $(D_2,\dashv_2, \vdash_2)$ be abelian dimonoids.  Dimonoids $(D_1,\dashv_1, \vdash_1)$ and $(D_2,\dashv_2, \vdash_2)$ are isomorphic if and only if semigroups $(D_1,\dashv_1)$ and $(D_2,\dashv_2)$ are isomorphic.
\end{proposition}

\begin{proof} It is immediate to show that if dimonoids  $(D_1,\dashv_1, \vdash_1)$ and $(D_2,\dashv_2, \vdash_2)$ are isomorphic, then semigroups $(D_1,\dashv_1)$ and $(D_2,\dashv_2)$ are isomorphic as well.
Conversely, let $\psi: D_1 \to D_2$ be an isomorphism from a semigroup $(D_1,\dashv_1)$ to a semigroup $(D_2,\dashv_2)$. Since $\psi(x\vdash_1 y) = \psi(y\dashv_1 x) = \psi(y)\dashv_2 \psi(x) = \psi(x)\vdash_2 \psi(y)$ for all $x,y\in D_1$, it follows that $\psi: D_1 \to D_2$ is an isomorphism from a dimonoid $(D_1,\dashv_1, \vdash_1)$ to a dimonoid $(D_2,\dashv_2, \vdash_2)$.
\end{proof}

\section{Two-element dimonoids and their automorphism groups}\label{2eldim}

In this section we describe, up to isomorphism, all two-element dimonoids and their automorphism groups.

\begin{theorem}
Up to isomorphism, there exist $8$ two-element dimonoids among which $3$ dimonoids are commutative. Also, up to isomorphism, there are $4$ abelian dimonoids of order $2$, and nonabelian dimonoids are divided into $2$ pairs of dual dimonoids. There exist exactly $5$ pairwise nonisomorphic two-element trivial dimonoids.
\end{theorem}

\begin{proof}

It is well-known that there are exactly five pairwise nonisomorphic semigroups having two elements: the multiplicative cyclic group $C_2=\{-1,1\}$, the linear semilattice $L_2=\{0,1\}$ with $\min$-operation, the null semigroup $O_2=\{0,1\}$ with zero $0$, the left zero semigroup $LO_2$ with operation $ab=a$, and the right zero semigroup $RO_2$ with operation $ab=b$.

In the sequel, we divide our investigation into cases. In the case of a semigroup $(S,*)$ we shall find all pairwise nonisomorphic dimonoids $(D,\dashv, \vdash)$ such that $(D,\dashv)$ is isomorphic to $(S,*)$.

\smallskip
{\noindent \bf Cases $C_2$ and $L_2$}.  According to Proposition~\ref{rl_id}, if a semigroup $(D,\dashv)$ possesses a left identity or a semigroup $(D,\vdash)$ possesses a right identity, then the operations of a dimonoid $(D,\dashv, \vdash)$ coincide. Therefore, up to isomorphism, there exist a unique dimonoid $(D,\dashv, \vdash)$ such that $(D,\dashv)\cong C_2$ or $(D,\vdash)\cong C_2$, and this dimonoid is the trivial dimonoid $C_2$. Similarly, $L_2$ is a unique dimonoid in the class of dimonoids $(D,\dashv, \vdash)$ such that $(D,\dashv)\cong L_2$ or $(D,\vdash)\cong L_2$. The trivial dimonoids $C_2$ and $L_2$ are commutative and abelian.

\smallskip
{\noindent \bf Case $LO_2$}. If $(D,\vdash)\cong RO_2$, then we obtain the abelian noncommutative dimonoid $LO_2\rbag RO_2$. According to Proposition~\ref{isoabelian}, $LO_2\rbag RO_2$ is a unique dimonoid in the class of abelian dimonoids $(D,\dashv, \vdash)$ such that  $(D,\dashv)\cong LO_2$ and $(D,\vdash)\cong RO_2$. It follows that $\Aut(LO_2\rbag RO_2)=\Aut(LO_2)=S_2\cong C_2$.

In the case $(D,\vdash)\cong O_2$,  we obtain the noncommutative nonabelian dimonoid $LO_2\rbag O_2$. By Proposition~\ref{isoLO}, $LO_2\rbag O_2$ is a unique dimonoid in the class of dimonoids $(D,\dashv, \vdash)$ such that $(D,\dashv)\cong LO_2$ and $(D,\vdash)\cong O_2$. It follows that $\Aut(LO_2\rbag O_2) = \Aut(O_2) = C_1$.

In the remaining case, we obtain the trivial noncommutative nonabelian dimonoid $LO_2$. By Proposition~\ref{isoLO}, $LO_2$ is a unique dimonoid in the class of dimonoids $(D,\dashv, \vdash)$ such that $(D,\dashv)\cong (D,\vdash)\cong LO_2$.

\smallskip
{\noindent \bf Case $RO_2$}. According to Proposition~\ref{rightcoin}, if $(D,\dashv)$ is a right zero semigroup, then the operations of a dimonoid $(D,\dashv, \vdash)$ coincide. Consequently, the trivial dimonoid $RO_2$ is a unique dimonoid in the class of  dimonoids $(D,\dashv, \vdash)$ such that $(D,\dashv)\cong RO_2$. This dimonoid is noncommutative and nonabelian, and it is dual to the dimonoid $LO_2$.

\smallskip
{\noindent \bf Case $O_2$}.  
If $(D,\vdash)$ is a right zero semigroup, then we obtain the noncommutative nonabelian dimonoid $O_2\rbag RO_2$, which is unique in the class of dimonoids $(D,\dashv, \vdash)$ such that $(D,\dashv)\cong O_2$ and $(D,\vdash)\cong RO_2$, in accordance with Proposition~\ref{isoRO}. This dimonoid is dual to the dimonoid $LO_2\rbag O_2$, and $\Aut(O_2\rbag RO_2) = C_1$.

According to Proposition~\ref{rightcoin}, if $(D,\vdash)$ is a left zero semigroup, then the operations of a dimonoid $(D,\dashv, \vdash)$ coincide. Therefore, there does not exist a dimonoid $(D,\dashv, \vdash)$ such that $(D,\dashv) \cong O_2$ and $(D,\vdash)\cong LO_2$.

In the final case, we obtain the trivial commutative abelian dimonoid $O_2$. It follows from Corollary~\ref{zerodim} that $O_2$ is a unique dimonoid in the class of dimonoids $(D,\dashv, \vdash)$ such that $(D,\dashv)\cong (D,\vdash)\cong O_2$.
\end{proof}

In the following table we present,  up to isomorphism, all two-element  dimonoids and their automorphism groups.

\begin{table}[ht]
    \centering

    \begin{tabular}{|c|c|c|c|c|c|c|c|c|}
        \hline
        $D$ & $C_2$  & $L_2$ & $O_2$  & $LO_2$ & $RO_2$ & $LO_2 \rbag RO_2$ & $LO_2 \rbag O_2$ & $O_2 \rbag RO_2$  \\
        \hline
        $\Aut(D)$ &  $C_1$ & $C_1$ & $C_1$ & $C_2$ & $C_2$ & $C_2$  & $C_1$ & $C_1$  \\
        \hline
    \end{tabular}

    \smallskip
    \caption{Two-element dimonoids and their automorphism groups}\label{tab:DopSG2}
\end{table}

\section{Three-element dimonoids and their automorphism groups}\label{3eldim}

In the remaining part of this paper, we focus on describing, up to isomorphism, all three-element dimonoids.

Among the $19683$ possible binary operations on a three-element set $S$, precisely $113$ are associative. In other words, there exist exactly $113$ distinct three-element semigroups. However, many of these semigroups are isomorphic, and as a result, there are essentially only $24$ pairwise nonisomorphic semigroups of order $3$, see \cite{Ch, G8, G9}.

Among these $24$ pairwise nonisomorphic semigroups of order $3$, there are $12$ commutative semigroups. The remaining $12$ pairwise nonisomorphic noncommutative semigroups are partitioned into pairs of dual semigroups. Moreover, the automorphism groups of dual semigroups coincide.

List of all pairwise nonisomorphic semigroups of order $3$ and their automorphism groups are presented in Table~\ref{tab:auts3} and Table~\ref{tab:autns3} taken from~\cite{G9}.

\begin{table}[ht]
    \centering

    \begin{tabular}{|c|c|c|c|c|c|c|c|c|c|c|c|c|}
        \hline
        $S$ & $C_3$ & $O_3$  & $\M_{2,2}$ & $C_2^{+1}$ &  $C_2^{\tilde{1}}$ & $\M_{3,1}$&  $O_2^{+1}$ & $O_2^{+0}$ & $L_3$ & $C_2^{+0}$    & $O_3^2$  & $O_3^1$ \\
        \hline
        $\Aut(S)$ &  $C_2$ & $C_2$ & $C_1$ & $C_1$ & $C_{1}$  & $C_{1}$ & $C_{1}$ & $C_{1}$ & $C_{1}$ & $C_{1}$ & $C_{2}$ & $C_{1}$  \\
        \hline
    \end{tabular}

    \smallskip
    \caption{Commutative semigroups of order  $3$ and their automorphism groups}\label{tab:auts3}
\end{table}

\begin{table}[ht]
    \centering
    \resizebox{14cm}{!}{
        \begin{tabular}{|c|c|c|c|c|c|c|}
            \hline
            $S$ & $LO_3$, $RO_3$  & $LO_2^{+0}$, $RO_2^{+0}$ &  $LO^{\sim 0}_{1\leftarrow2}$, $RO^{\sim 0}_{1\leftarrow2}$ & $LO_2^{+1}$, $RO_2^{+1}$ & $LOB_3$, $ROB_3$  & $LO_{2\leftarrow 3}$, $RO_{2\leftarrow 3}$ \\
            \hline
            $\Aut(S)$ & $S_3$ & $C_2$ & $C_1$ & $C_{2}$ & $C_{1}$ & $C_{2}$ \\
            \hline
        \end{tabular}
    }
    \smallskip
    \caption{Noncommutative three-element  semigroups and their automorphism groups}\label{tab:autns3}
\end{table}

\subsection{Commutative three-element dimonoids}

The classification of three-element commutative dimonoids will be based on our results concerning the classification of three-element doppelsemigroups from~\cite{GR1}.

A {\em doppelsemigroup} is an algebraic structure $(D,\dashv,\vdash)$ consisting of a nonempty set $D$ equipped with two associative binary operations $\dashv$ and $\vdash$ satisfying the axiom $(D_2)$ and the following axiom:
$$(x \dashv y) \vdash z = x \dashv (y \vdash z) \hspace{10mm} (D_4).$$

For a doppelsemigroup $(D,\dashv,\vdash)$, if $\mathbb S$ and $\mathbb T$ denote the semigroups $(D,\dashv)$ and $(D,\vdash)$, respectively, then $\mathbb S \between \mathbb T$ stands for the doppelsemigroup $(D,\dashv,\vdash)$. 

%Let $(D_1,\dashv_1, \vdash_1)$ be such a doppelsemigroup that for each doppelsemigroup $(D_2,\dashv_2, \vdash_2)$ the isomorphisms $(D_2,\dashv_2)\cong (D_1,\dashv_1)$ and $(D_2,\vdash_2)\cong (D_1,\vdash_1)$ imply $(D_2,\dashv_2, \vdash_2)\cong(D_1,\dashv_1, \vdash_1)$. If $\mathbb S$ and $\mathbb T$ are model semigroups of classes of semigroups isomorphic to
%$(D_1,\dashv_1)$ and $(D_1, \vdash_1)$, respectively, then by $\mathbb S\between\mathbb T$  we denote a model doppelsemigroup of the class of doppelsemigroups isomorphic to $(D_1,\dashv_1, \vdash_1)$.

In \cite{GR1}, the problem of classifying all doppelsemigroups with at most three elements up to isomorphism was completely solved. According to Proposition~1 from \cite{ZhDADM2016}, every commutative dimonoid is a doppelsemigroup. Therefore, in order to describe all three-element commutative dimonoids up to isomorphism, it suffices to select those dimonoids among the commutative pairwise nonisomorphic doppelsemigroups of order $3$.

The following Table~\ref{tab:commdop3} of all pairwise nonisomorphic nontrivial commutative three-element doppelsemigroups and their automorphism groups is taken from~\cite{GR1}.

\begin{table}[ht]
    \centering
    \resizebox{13cm}{!}{
        \begin{tabular}{|c|c|c|c|c|c|c|}
             \hline
            $D$ & $C_3\between C_3^{-1}$  &  $O_3\between \M_{3,1}$ & $O_3\between O_2^{+1}$ & $O_3\between O_2^{+0}$  & $O_3\between L_3$ & $O_3\between C_2^{+0}$  \\
              \hline
            $\Aut(D)$ & $C_1$  & $C_1$ & $C_{1}$ & $C_{1}$ & $C_{1}$ & $C_1$ \\
            \hline
            \hline
            $D$   & $O_3\between O_3^2$  & $O_3\between O_3^1$ &   $\M_{2,2}\between C_2^{+1}$  &  $\M_{2,2}\between C_2^{\tilde{1}}$ & $C_2^{+1}\between C_2^{\tilde{1}}$ & $C_2^{+1}\between \M_{2,2}$ \\
            \hline
            $\Aut(D)$  & $C_2$ & $C_1$ & $C_{1}$ & $C_{1}$ & $C_{1}$ & $C_{1}$\\
            \hline
            \hline
            $D$ & $C_2^{\tilde{1}}\between \M_{2,2} $ & $C_2^{\tilde{1}}\between C_2^{+1} $ & $\M_{3,1}\between O_2^{+1}$ & $\M_{3,1}\between O_3$  & $O_2^{+1}\between \M_{3,1}$ & $O_2^{+1}\between O_3$  \\
\hline
 $\Aut(D)$ & $C_1$ & $C_1$ & $C_1$ & $C_{1}$ & $C_{1}$ & $C_{1}$\\
 \hline
\hline
 $D$ & $(O_2\between L_2)^{+0}$  &  $ O_2^{+0}\between O_3$ & $L_3 \between O_3$ & 
   $(L_2\between O_2)^{+0}$ & $(C_2\between C_2^{-1})^{+0}$ & $C_2^{+0}\between O_3$     \\
             \hline
            $\Aut(D)$ & $C_1$ & $C_1$ & $C_1$ & $C_{1}$ & $C_{1}$ & $C_{1}$\\
            \hline
            \hline
            $D$ & $O_3^2\between O_3^1$ & $O_3^2\between O_3$  & $O_3^a\between O_3^b$ & $O_3^1\between O_3^2$ &  $O_3^1\between O_3$   \\
             \cline{1-6}
            $\Aut(D)$ & $C_1$ & $C_2$ & $C_1$ & $C_{1}$ & $C_{1}$  \\
            \cline{1-6}
        \end{tabular}

    }
    \smallskip
    \caption{Three-element nontrivial commutative doppelsemigroups and their automorphism groups}\label{tab:commdop3}
\end{table}

We begin by establishing several auxiliary propositions.

\begin{proposition}\label{lnulldim} Let $(D,\dashv, \vdash)$ be a doppelsemigroup such that $(D,\dashv)$ is a null semigroup with zero $0$. A doppelsemigroup $(D,\dashv, \vdash)$ is a dimonoid if and only if $D \vdash D \vdash D = \{0\}$.
\end{proposition} 
\begin{proof}
Taking into account that for a doppelsemigroup $(D,\dashv, \vdash)$ an element $0\in D$ is a zero of a semigroup $(D,\dashv)$ if and only if $0$ is a zero of a semigroup $(D,\vdash)$, see \cite{GR1}, we conclude that the axioms $(D_1)$ and $(D_2)$ of a dimonoid hold:  
$$(x \dashv y) \dashv z = 0 = x \dashv (y \vdash z),\ \ \ (x \vdash y) \dashv z = 0 = x \vdash 0 = x \vdash (y \dashv z).$$
Since $(x \dashv y) \vdash z = 0 \vdash z = 0$ for any $x,y,z\in D$, we conclude that the axiom $(D_3)$ of a dimonoid holds if and only if $x \vdash (y \vdash z)=0$ for any $x,y,z\in D$, that is $D \vdash D \vdash D = \{0\}$.  
\end{proof}

Dually, one can prove the following proposition.

\begin{proposition}\label{rnulldim} Let $(D,\dashv, \vdash)$ be a  doppelsemigroup such that $(D,\vdash)$ is a null semigroup with zero $0$. A doppelsemigroup $(D,\dashv, \vdash)$ is a  dimonoid if and only if $D \dashv D \dashv D = \{0\}$.
\end{proposition}

The following theorem provides a complete classification of all pairwise nonisomorphic commutative dimonoids of order~$3$.

\begin{theorem}
Up to isomorphism, there exist $14$  three-element commutative dimonoids among which $12$ trivial dimonoids and a pair of nontrivial nonabelian dual dimonoids. 
\end{theorem}

\begin{proof}

Since a trivial dimonoid $(D,\dashv, \dashv)$ is commutative if and only if a semigroup $(D,\dashv)$ is commutative, we obtain that, up to isomorphism, there exist $12$ trivial commutative dimonoids, see~Table~\ref{tab:auts3}. 

\smallskip

Our further investigation is carried out by distinguishing several cases.

\smallskip

{\noindent \bf Case 1}. Consider the doppelsemigroups $O_3\between \M_{3,1}$ and $\M_{3,1}\between O_3$. Recall that
$\M_{3,1}=\{a,a^2,a^3\ |\ a^4=a^3\}$ is a monogenic semigroup of index $3$ and period $1$ with zero $a^3$. Since $\M_{3,1}*\M_{3,1}*\M_{3,1}=\{a^3\}$,  we conclude that $O_3\between \M_{3,1}$ and $\M_{3,1}\between O_3$ are (nonabelian dual) dimonoids according to Propositions~\ref{lnulldim} and~\ref{rnulldim}.
These dimonoids are examples of commutative nonabelian dimonoids.

\smallskip

{\noindent \bf Case 2}.   According to Proposition~\ref{rl_id} for a dimonoid $(D,\dashv, \vdash)$, if a semigroup $(D,\dashv)$ contains a left identity or a semigroup $(D,\vdash)$ contains a right identity, then the operations of a dimonoid $(D,\dashv, \vdash)$ coincide. Therefore, the doppelsemigroups $C_3\between C_3^{-1}$, $O_3\between O_2^{+1}$, $O_3\between L_3$, $O_3\between C_2^{+0}$, $\M_{2,2}\between C_2^{+1}$, $C_2^{+1}\between C_2^{\tilde{1}}$, $C_2^{+1}\between \M_{2,2}$, $C_2^{\tilde{1}}\between C_2^{+1} $, $\M_{3,1}\between O_2^{+1}$, $O_2^{+1}\between \M_{3,1}$, $O_2^{+1}\between O_3$, $(O_2\between L_2)^{+0}$,  $L_3 \between O_3$, $(L_2\between O_2)^{+0}$, $(C_2\between C_2^{-1})^{+0}$, and $C_2^{+0}\between O_3$ can not be dimonoids.

\smallskip

{\noindent \bf Case 3}. Consider the doppelsemigroups $O_3\between O_2^{+0}$ and $O_2^{+0}\between O_3$. Let $0$ and $z$ be  zeros of the semigroup $O_3$ and $O_2$, respectively. Taking into account that $O_2^{+0}*O_2^{+0}*O_2^{+0} = \{0,z\}\ne \{0\}$,  we conclude according to Propositions~\ref{lnulldim} and~\ref{rnulldim} that $O_3\between O_2^{+0}$ and $O_2^{+0}\between O_3$ are not dimonoids. 

\smallskip

{\noindent \bf Case 4}. Consider the doppelsemigroups $O_3\between O_3^2$ and $O_3^2\between O_3$. Recall that $O_3^2$ is a nonlinear semilattice isomorphic to the semigroup $\{a,b, 0\}$ with the operation $*$:
$$x* y=\begin{cases}
x,\text{ if } y=x\in\{a,b\} \\
0,\ \text{otherwise.}
\end{cases}$$
Since $O_3^2*O_3^2*O_3^2 = O_3^2\ne \{0\}$,  we conclude  that $O_3\between O_3^2$ and $O_3^2\between O_3$ are not dimonoids by Propositions~\ref{lnulldim} and~\ref{rnulldim}.

\smallskip

{\noindent \bf Case 5}. Consider the doppelsemigroups $O_3\between O_3^1$ and $O_3^1\between O_3$. Recall that $O_3^1$ is isomorphic to the semigroup $\{a, b, 0\}$ with the operation $*$:
$$x* y=\begin{cases}
x,\text{ if } y=x=a \\
0,\ \text{otherwise.}
\end{cases}$$
Since $O_3^1*O_3^1*O_3^1 = \{0,a\}\ne \{0\}$,  we conclude  that $O_3\between O_3^1$ and $O_3^1\between O_3$ are not dimonoids according to Propositions~\ref{lnulldim} and~\ref{rnulldim}. 

\smallskip

%A semigroup $(S,*)$ is called {\em globally idempotent} if $S*S = S$.

{\noindent \bf Case 6}. Consider the doppelsemigroups $O_3^2\between O_3^1$ and $O_3^1\between O_3^2$. 
According to Lemma~3 of~\cite{ZhAL2011} for a dimonoid $(D,\dashv, \vdash)$, if  $(D,\dashv)$ is a semilattice, then the operations of a dimonoid $(D,\dashv, \vdash)$ coincide. Since $O_3^2$ is a semilattice, the doppelsemigroup $O_3^2\between O_3^1$ can not be a dimonoid. The doppelsemigroup $O_3^1\between O_3^2$ also cannot be a dimonoid either, because otherwise $O_3^2\between O_3^1$ would be its dual dimonoid.

\smallskip

{\noindent \bf Case 7}. Consider the doppelsemigroups $\M_{2,2}\between C_2^{\tilde{1}}$ and $C_2^{\tilde{1}} \between  \M_{2,2}$. Recall that $\M_{2,2}=\{a,a^2,a^3\ |\ a^4=a^2\}$ with operation $\dashv$ is a monogenic semigroup of index 2 and period 2 and $C_2^{\tilde{1}} =\{a^2, a^3\}^{\tilde{1}}$ with operation $\vdash$ is a semigroup obtained from the cyclic group $\{a^2,a^3\}$ with identity $a^2$ by adjoining an element $a$ with $a\vdash s = s\vdash a = s$  for $s\in \{a^2,a^3\}$ and $a\vdash a=a^2$. Taking into account that
$(a\dashv a)\dashv a^2 = a^2 \dashv a^2 = a^4 = a^2$ and $a\dashv (a\vdash a^2) = a \dashv a^2 = a^3 \ne a^2$,
we conclude that $\M_{2,2}\between C_2^{\tilde{1}}$ is not a dimonoid. By analogy $C_2^{\tilde{1}} \between  \M_{2,2}$ is not a dimonoid.

\smallskip

{\noindent \bf Case 8}. Consider the last doppelsemigroup $O_3^a\between O_3^b$. Recall that $O_3^a\between O_3^b$ is  the doppelsemigroup $(\{a, b, 0\}, *_a, *_b)$, where for $t\in \{a,b\}$
$$x*_t y=\begin{cases}
x,\text{ if } y=x=t \\
0,\ \text{otherwise.}
\end{cases}$$
Taking into account that
$(a *_a a)*_a a = a *_a a = a$ and $a *_a  (a *_b a ) = a *_a 0 = 0 \ne a$,
we conclude that $O_3^a\between O_3^b$ is not a dimonoid. 
\end{proof}

\bigskip

In the following Table~\ref{tab:commdim3} we present,  up to isomorphism, all three-element  commutative dimonoids and their automorphism groups.

\begin{table}[ht]
\centering
\resizebox{15cm}{!}{
    \begin{tabular}{|c|c|c|c|c|c|c|c|c|c|c|c|c|c|c|}
        \hline
        $D$ & $C_3$ & $O_3$  & $\M_{2,2}$ & $C_2^{+1}$ &  $C_2^{\tilde{1}}$ & $\M_{3,1}$&  $O_2^{+1}$ & $O_2^{+0}$ & $L_3$ & $C_2^{+0}$    & $O_3^2$  & $O_3^1$ & $\M_{3,1} \rbag O_3$ & $O_3 \rbag \M_{3,1}$\\
        \hline
        $\Aut(D)$ &  $C_2$ & $C_2$ & $C_1$ & $C_1$ & $C_{1}$  & $C_{1}$ & $C_{1}$ & $C_{1}$ & $C_{1}$ & $C_{1}$ & $C_{2}$ & $C_{1}$ & $C_{1}$ & $C_{1}$ \\
        \hline
    \end{tabular}
}
\smallskip
\caption{Commutative three-element dimonoids and their automorphism groups}\label{tab:commdim3}
\end{table}

\subsection{Abelian  three-element dimonoids}

The following theorem provides a complete classification of all pairwise nonisomorphic abelian dimonoids of order~$3$.

Recall that a semigroup $(S, *)$ is called {\em right}  {\em commutative}~\cite{ZhADM2013}, if it satisfies the identity $s*x*y = s*y*x$  for all $s, x, y \in S$.

\begin{theorem}
Up to isomorphism, there exist $17$  three-element abelian dimonoids among which $12$ commutative trivial dimonoids and $5$ noncommutative nontrivial dimonoids.
\end{theorem}

\begin{proof}
Since a trivial dimonoid $(D,\dashv, \dashv)$ is abelian if and only if a semigroup $(D,\dashv)$ is commutative, we obtain that, up to isomorphism, there exist $12$ trivial abelian dimonoids, see Table~\ref{tab:auts3}. 

Let $(D, \dashv)$ be an arbitrary semigroup and $(D, \vdash)$ be a dual semigroup to $(D, \dashv)$. According to Lemma 3 of~\cite{ZhADM2013}, an algebraic structure $(D, \dashv, \vdash)$ is an abelian dimonoid if and only if $(D, \dashv)$ is a right commutative semigroup.
If an abelian dimonoid $(D, \dashv, \vdash)$ has a commutative semigroup $(D, \dashv)$, then for all $x, y \in D$, it holds that $x \vdash y = y \dashv x = x \dashv y$. In this case, both operations coincide, and the dimonoid is trivial.

Consider two abelian dimonoids $(D_1, \dashv_1, \vdash_1)$ and $(D_2, \dashv_2, \vdash_2)$.  By Proposition~\ref{isoabelian}, dimonoids $(D_1,\dashv_1, \vdash_1)$ and $(D_2,\dashv_2, \vdash_2)$ are isomorphic if and only if semigroups $(D_1,\dashv_1)$ and $(D_2,\dashv_2)$ are isomorphic.

From the previous considerations it follows that  the task of describing all pairwise nonisomorphic nontrivial abelian three-element dimonoids reduces to the task of recognizing right commutative semigroups among the nontrivial noncommutative semigroups listed in Table~\ref{tab:autns3}.

It was proved in~\cite{Gdim1} that the semigroups $LO_3$, $LO_{2\leftarrow 3}$, $LOB_3$,  $LO^{\sim 0}_{1\leftarrow2}$, and $LO_2^{+0}$ are right commutative. 

It follows directly from the definition of a right commutative semigroup that a noncommutative semigroup containing a left identity cannot be right commutative. Therefore, the semigroups $RO_3$, $RO_2^{+0}$, $RO^{\sim 0}_{1\leftarrow2}$, $LO_2^{+1}$, $RO_2^{+1}$, and $ROB_3$ are not right commutative. Consider the remaining semigroup $RO_{2 \leftarrow 3}$, which contains two right zeros. Denote these zeros by $a$ and $b$. For any $s \in RO_{2 \leftarrow 3}$, it holds that $sab = b \neq a = sba$, and therefore, the semigroup $RO_{2 \leftarrow 3}$ is not right commutative.

We conclude that up to isomorphism there exist $5$ abelian noncommutative nontrivial dimonoids: 
$LO_3\rbag RO_3$, $LO_{2\leftarrow 3}\rbag RO_{2\leftarrow 3}$, $LOB_3\rbag ROB_3$,   $LO^{\sim 0}_{1\leftarrow 2}\rbag RO^{\sim 0}_{1\leftarrow 2}$, and $(LO_2\rbag RO_2)^{+0}= LO_2^{+0}\rbag RO_2^{+0}$. 
\end{proof}

Based on the results of~\cite{Gdim1} concerning the automorphism groups of abelian noncommutative dimonoids, Table~\ref{tab:noncommabdim3} lists, up to isomorphism, all  abelian noncommutative nontrivial three-element dimonoids and their corresponding automorphism groups.

\begin{table}[ht]
    \centering
    \resizebox{12.5cm}{!}{
        \begin{tabular}{|c|c|c|c|c|c|c|}
            \hline
            $D$  & $LO_3\rbag RO_3$  & $LO_{2\leftarrow 3}\rbag RO_{2\leftarrow 3}$ & $LOB_3\rbag ROB_3$ &   $LO^{\sim 0}_{1\leftarrow 2}\rbag RO^{\sim 0}_{1\leftarrow 2}$ & $(LO_2\rbag RO_2)^{+0}$  \\
             \cline{1-6}
            $\Aut(D)$ & $S_3$ & $C_2$ & $C_1$ & $C_{1}$ & $C_{2}$  \\
            \cline{1-6}
        \end{tabular}

    }
    \smallskip
    \caption{ Abelian noncommutative  nontrivial $3$-element dimonoids and their automorphism groups}\label{tab:noncommabdim3}
\end{table}

\subsection{Nonabelian noncommutative three-element dimonoids}

Based on the results of~\cite{Gdim1} concerning noncommutative nonabelian dimonoids and their automorphism groups and properties of dual dimonoids,
 Table~\ref{tab:noncommdim3} lists  some pairwise nonisomorphic   noncommutative nonabelian nontrivial three-element dimonoids and their automorphism groups.

\begin{table}[ht]
    \centering
    \resizebox{15cm}{!}{
        \begin{tabular}{|c|c|c|c|c|c|c|c|}
             \hline
            $D$ & $LO_3\rbag O_3$  &  $LO_{2\leftarrow 3}\rbag O_3$ & $LO_3\rbag RO_{2\leftarrow 3}$ & $LO_3\rbag LO_{2\leftarrow 3}$ &  $LOB_3\rbag O_3^1$  & $LO^{\sim 0}_{1\leftarrow2}\rbag O_3^1$   &  $(LO_2\rbag O_2)^{+0}$  \\
              \hline
            $\Aut(D)$ & $C_2$  & $C_1$ & $C_{2}$ & $C_{2}$ & $C_{1}$ & $C_{1}$ & $C_1$ \\
            \hline
            \hline
            $D$   & $O_3\rbag RO_3$   & $O_3\rbag RO_{2\leftarrow 3}$ &  $LO_{2\leftarrow 3}\rbag RO_3$ &  $RO_{2\leftarrow 3}\rbag RO_3$  & $O_3^1 \rbag ROB_3$ & $O_3^1 \rbag RO^{\sim 0}_{1\leftarrow2}$  & $(O_2\rbag RO_2)^{+0}$ \\
            \hline
            $\Aut(D)$  & $C_2$ & $C_1$ & $C_{2}$ & $C_{2}$ & $C_{1}$ & $C_{1}$ & $C_{1}$\\
            \hline
           
        \end{tabular}

    }
    \smallskip
    \caption{ Nonabelian  noncommutative nontrivial $3$-element  dimonoids and their automorphism groups}\label{tab:noncommdim3}
\end{table}

Since a trivial dimonoid $(D,\dashv, \dashv)$ is nonabelian if and only if a semigroup $(D,\dashv)$ is noncommutative, we obtain that, up to isomorphism, there exist $12$ trivial nonabelian noncommutative dimonoids, see Table~\ref{tab:autns3}. 
It follows that we have proved the following theorem.

\begin{theorem}
There exist at least $26$ pairwise nonisomorphic nonabelian noncommutative three-element dimonoids among which there are exactly $6$ pairs of trivial dual dimonoids and at least $7$ pairs of nontrivial dual dimonoids.
\end{theorem}

Solving the following problem will provide a complete classification, up to isomorphism, of all three element dimonoids.

\begin{problem}\label{3noncomdim} Give a complete classification, up to isomorphism, of all noncommutative nonabelian nontrivial three-element dimonoids.
%Is each such dimonoid isomorphic to one from Table~\ref{tab:noncommdim3}?
\end{problem}

\end{document}